\documentclass[11pt]{amsart}

\usepackage{amsfonts, amssymb, amscd}
\usepackage{amsmath}
\usepackage{verbatim}
\usepackage{amssymb}
\usepackage{mathrsfs}
\usepackage{graphicx}
\usepackage{cite}
\usepackage[all]{xy}
\usepackage{tikz}
\usepackage[toc,page]{appendix}
\usepackage{tikz-cd}

\usepackage{graphicx}
\usepackage{float}

\newtheorem{theorem}{Theorem}[section]
\newtheorem{lemma}[theorem]{Lemma}
\newtheorem{proposition}[theorem]{Proposition}
\newtheorem{definition}[theorem]{Definition}
\newtheorem{example}[theorem]{Example}
\newtheorem{corollary}[theorem]{Corollary}
\newtheorem{remark}[theorem]{Remark}

\numberwithin{equation}{section}
\begin{document}
\title{On $\mathrm{H}-$trivial line bundles on toric DM stacks of dim $\geq3$}

\begin{abstract}
We study line bundles on smooth toric Deligne-Mumford stacks $\mathbb{P}_{\mathbf{\Sigma}}$ of arbitrary dimension. We give a sufficient condition for when infinitely many line bundles on $\mathbb{P}_{\mathbf{\Sigma}}$ have trivial cohomology. In dimension three, this sufficient condition is also a necessary condition under the technical assumption that $\mathbf{\Sigma}$ has no more than one pair of collinear rays.
\end{abstract}

\author{Lev Borisov}
\address{Department of Mathematics\\
Rutgers University\\
Piscataway, NJ 08854} \email{borisov@math.rutgers.edu}

\author{Chengxi Wang}
\address{Department of Mathematics\\
University of California Los Angeles\\
Los Angeles, CA 90095} \email{chwang@math.ucla.edu}

\maketitle

\tableofcontents
\section{Introduction}\label{section1}
Derived categories of coherent sheaves on algebraic varieties and stacks have attracted significant attention over the last several decades. Some of these categories can be generated by so-called full exceptional collections, often built from line bundles,
\cite{CM,HP1,HP2,I,Ka,Ki,M}.

\smallskip
A defining feature of an exceptional collection $(\mathcal L_1,\ldots, \mathcal L_n)$
of line bundles is that for $i<j$ the cohomology spaces $H^k(\mathcal L_j \otimes \mathcal L_i^{-1})$ vanish for all $k$. This prompts a problem, which is rather intriguing, even outside of the context of derived categories:

\smallskip
{\em For a given variety/stack, characterize all line bundles on it with trivial cohomology.}

\smallskip
These line bundles are called $\mathrm{H}-$trivial in \cite{W} and the current paper and \emph{immaculate} in \cite{Aetal}.
Unfortunately, the problem of  classification of all $H$-trivial line bundles is too hard  in many cases. However, the following, more tractable question is the first step towards this goal.

\smallskip
{\em Does a given  variety/stack admit infinitely many  $\mathrm{H}-$trivial line bundles?}

\smallskip
Our main focus in this paper is on smooth proper toric Deligne-Mumford (DM) stacks, where one can hope to translate the answer to the above question into a combinatorial criterion. For example, the
paper \cite{W} gives such criterion in the dimension two case.

\begin{theorem}\label{2}
\cite{W} Let $\mathbb{P}_{\mathbf{\Sigma}}$ be a proper smooth dimension two toric Deligne-Mumford stack associated to a complete stacky fan $\mathbf{\Sigma}=(\Sigma, \{v_i\}_{i=1}^{n})$. Then
there are infinitely many $\mathrm{H}-$trivial line bundles on $\mathbb{P}_{\mathbf{\Sigma}}$ if and only if there exists $\{i, j\}\subset\{1,2,\ldots,n\}$ such that $v_i$ and $v_j$ are collinear.
\end{theorem}

\begin{remark}
The importance of pairs of collinear rays was already observed in the paper \cite{P} in the related context of vanishing of cohomology of divisorial sheaves on toric varieties.
\end{remark}

This paper came out of a project to extend the result of \cite{W} to toric DM stacks  $\mathbb{P}_{\mathbf{\Sigma}}$ in higher dimensions.
While we have made significant progress, a complete criterion is still elusive, even in dimension three.

\smallskip
Our first result gives a sufficient condition for infinitude of $\mathrm{H}-$trivial line bundles.  We associate with each $\mathbf{\Sigma}-$piecewise linear function $\psi$ a convex polytope $\Lambda_{\psi}$ in the real span of lattice of characters by looking at the corresponding restrictions to the maximum cones of $\mathbf{\Sigma}$, see Definition \ref{deL}.

\smallskip
\smallskip\noindent
{\bf Theorem \ref{if}.}
Let $\mathbb{P}_{\mathbf{\Sigma}}$ be a proper smooth dimension $m$ toric DM stack associated to a complete stacky fan $\mathbf{\Sigma}=(\Sigma, \{v_i\}_{i=1}^{n})$. If there exists a $\mathbf{\Sigma}-$piecewise linear function $\psi$ such that $0<\mathrm{dim}(\Lambda_{\psi})<m$, then
there are infinitely many $\mathrm{H}-$trivial line bundles on $\mathbb{P}_{\mathbf{\Sigma}}$.

\begin{remark}
If $\,\mathbb{P}_{\mathbf{\Sigma}}$ admits a non-trivial fibration structure
$
\mathbb{P}_{\mathbf{\Sigma}}\to  \mathbb{P}_{\mathbf{\Sigma_1}}
$
then a pull-back of a non-trivial line bundle from $  \mathbb{P}_{\mathbf{\Sigma_1}}$ provides the requisite  $\mathbf{\Sigma}-$piecewise linear function $\psi$ of Theorem \ref{if}. However, such  $\psi$
may exist without a fibration structure, so  existence of $\psi$ can be viewed as some weaker version of a fibration structure.
\end{remark}

\smallskip

We also get a partial converse of Theorem \ref{if}  for smooth toric varieties and DM stacks in dimension three.

\smallskip

{\bf Theorem \ref{oneonly}.}
Let $\mathbb{P}_{\mathbf{\Sigma}}$ be a proper smooth dimension three toric DM stack associated to a complete stacky fan $\mathbf{\Sigma}=(\Sigma, \{v_i\}_{i=1}^{n})$. Assume there exists no more than one pair of collinear rays in $\mathbf{\Sigma}$. Then
there are infinitely many $\mathrm{H}-$trivial line bundles on $\mathbb{P}_{\mathbf{\Sigma}}$ if and only if there exists a $\mathbf{\Sigma}-$piecewise linear function $\psi$ such that $0<\mathrm{dim}(\Lambda_{\psi})<3$.

\smallskip
\smallskip

The paper is organized as follows.
In Section \ref{stack}, we give an overview of smooth toric DM stacks, their Picard groups and the cohomology of line bundles on the stacks. Then we define forbidden cones and forbidden sets and state the first main result Theorem \ref{if}. Section \ref{H} focuses on the proof of Theorem \ref{if}. We first exhibit an important way of producing infinitely many $\mathrm{H}-$trivial line bundles in Proposition \ref{if1}. Then we relate it to the existence of a $\mathbf{\Sigma}-$piecewise linear function $\psi$ such that $0<\mathrm{dim}(\Lambda_{\psi})<m$. In Section \ref{three}, we consider the case of dimension three. We prove a sufficient and necessary condition for the existence of infinitely many $\mathrm{H}-$trivial line bundles under the assumption that there is no more than one pair of collinear rays in $\mathbf{\Sigma}$. Section \ref{comments} describes our current state of knowledge and states some open questions related to infinitude of $\mathrm{H}-$trivial line bundles on toric DM stacks.
\smallskip

\noindent{\it Acknowledgements.}
L. Borisov was partially supported by
NSF grant DMS-1601907. We thank Markus Perling for pointing out the paper \cite{P} to us.
\section{Line bundles on toric DM stacks and their cohomology}\label{stack}
In this section, we introduce toric DM stacks $\mathbb{P}_{\mathbf{\Sigma}}$ and their Picard groups $\mathrm{Pic} (\mathbb{P}_{\mathbf{\Sigma}})$, and describe the cohomology  spaces of line bundles on $\mathbb{P}_{\mathbf{\Sigma}}$. These results are well known but we need to state them to set up the notation and terminology and to help the reader who is not familiar with previous work in the area.

\smallskip
To avoid the technicalities of the derived Gale duality of \cite{BCS}, we consider a lattice $N$ which is a free abelian group of finite rank. Let $\Sigma$ be a complete simplicial fan in $N$. We choose a lattice point $v$ in each of the one-dimensional cones of $\Sigma$. If $\Sigma$ has $n$ one-dimensional cones, we get a complete stacky fan $\mathbf{\Sigma}=(\Sigma, \{v_i\}_{i=1}^{n})$, see \cite{BCS}.
The toric DM stack $\mathbb{P}_{\mathbf{\Sigma}}$ associated to this stacky fan $\mathbf{\Sigma}$ is constructed in \cite{BCS} as a stack version of the homogeneous coordinate ring construction of \cite{C}. The description of line bundles on the DM stacks is analogous to the description of the Picard group that was given in \cite{D,F}.
By \cite{V}, we know the line bundles on $\mathbb{P}_{\mathbf{\Sigma}}$ are in bijection with collections of integers, up to global linear functions, as described below.

\begin{proposition}\label{pic}
The Picard group of $\mathbb{P}_{\mathbf{\Sigma}}$ is isomorphic to the quotient of $\mathbb{Z}^n$ with basis $\{E_i\}_{i=1}^{n}$ by the subgroup of elements of the form $\sum_{i=1}^{n}(w\cdot v_i)E_i$ for all $w$ in the character lattice $M=N^*$. We use  the notation $\mathcal{O}(\sum_{i=1}^{n}r_iE_i)$ to denote these line bundles.
\end{proposition}
\begin{proof}
See Proposition 3.3 in \cite{K}.
\end{proof}

Now we remind the reader how to calculate the cohomology of a line bundle $\mathcal{L}$ on $\mathbb{P}_{\mathbf{\Sigma}}$. For each $\mathbf{r}=(r_i)_{i=1}^{n}\in \mathbb{Z}^n$, we define $Supp(\mathbf{r})$ to be the abstract simplicial complex on $n$ vertices $\{1,\ldots,n\}$ as follows
\begin{equation*}
\begin{split}
Supp(\mathbf{r})=&\{J \subseteq \{1,\ldots,n\}| r_i\geq 0 \text{ for all }i\in J\\
&\text{ and there exists a cone of $\Sigma$ containing all }v_i, i\in J \}.
\end{split}
\end{equation*}
The following proposition gives a description of the cohomology of a line bundle $\mathcal{L}$ on $\mathbb{P}_{\mathbf{\Sigma}}$ in terms of the reduced simplicial homology spaces of $Supp(\mathbf{r})$.
\begin{proposition}\label{Coho}
\cite{K} Let $\mathcal{L} \in \mathrm{Pic}(\mathbb{P}_{\mathbf{\Sigma}})$. Then $$\mathrm{H}^j(\mathbb{P}_{\mathbf{\Sigma}},\mathcal{L})=\bigoplus \mathrm{H}^{red}_{rk(N)-j-1}(Supp(\mathbf{r})),$$ where the sum is over all $\mathbf{r}=(r_i)_{i=1}^{n}\in \mathbb{Z}^n$ such that $\mathcal{O}(\sum_{i=1}^{n}r_iE_i)\cong \mathcal{L}$.
\end{proposition}
\begin{proof}
See Proposition 4.1 in \cite{K}.
\end{proof}
\begin{remark}\label{extreme}
We have $\mathrm{H}^0(\mathcal{L})\neq0$ if and only if there exists $\mathbf{r}\in \mathbb{Z}_{\geq 0}^n$ such that $\mathcal{O}(\sum_{i=1}^{n}r_iE_i)\cong \mathcal{L}$. The other extreme case of $\mathrm{H}^{rk(N)}(\mathcal{L})\neq 0$  happens when the simplicial complex $Supp(\mathbf{r})=\{\emptyset\}$, i.e. when there exists $\mathbf r$ such that $\mathcal{O}(\sum_{i=1}^{n}r_iE_i)\cong \mathcal{L}$ with all $r_i\leq -1$.
\end{remark}

\begin{remark}\label{coho2}
Let $\mathcal{L}\cong \mathcal{O}(\sum_{i=1}^{n}a_iE_i) $ be a line bundle in $\mathrm{Pic} (\mathbb{P}_{\mathbf{\Sigma}})$. Assume there is another expression $\mathcal{L} \cong \mathcal{O}(\sum_{i=1}^{n}r_iE_i)$. Then by Proposition \ref{pic}, there exists an element $f\in M$ such that $r_i=a_i+f \cdot v_i$ for $i=1,\ldots,n$. Thus the cohomology of $\mathcal{L}$ can also be written as
$$\mathrm{H}^j(\mathbb{P}_{\mathbf{\Sigma}},\mathcal{L})=\bigoplus_{f\in N^*}\mathrm{H}^{red}_{rk(N)-j-1}(Supp(\mathbf{r}_f)),$$where $\mathbf{r}_f=(a_i+f\cdot v_i)_{i=1}^{n}$.
\end{remark}

In this paper, our primary objects of interest are $\mathrm{H}-$trivial line bundles which we define below.

\begin{definition}\label{defH}
Let $\mathcal{L} $ be a line bundle in $\mathrm{Pic} (\mathbb{P}_{\mathbf{\Sigma}})$. We say that $\mathcal{L}$ is $\mathrm{H}-$trivial iff $\mathrm{H}^j(\mathbb{P}_{\mathbf{\Sigma}},\mathcal{L})=0$ for all $j\geq 0$.
\end{definition}

A combinatorial criterion for $\mathrm{H}-$triviality is given in terms of forbidden sets introduced below, see \cite{K,E}.

\begin{definition}\label{deforbiddenset}
For every subset $I \subseteq \{1,\ldots,n\}$, we denote $C_I$ to be the simplicial complex $Supp(\mathbf{r})$ where $r_i=-1$ for $i \notin I$ and $r_i=0$ for $i\in I$. Let $\Delta=\{I \subseteq \{1,\ldots,n\}| C_I \text{ has nontrivial reduced homology}\}$.  For each $I\in \Delta$, the forbidden set associated to $I$ is defined by $$FS_{I}:=\{\mathcal{O}(\sum_{i\notin I}(-1-r_i)E_i+\sum_{i\in I}r_iE_i)|r_i\in \mathbb{Z}_{\geq 0} \text{ for all } i\}.$$
\end{definition}

\begin{remark}
Since $\mathbf{\Sigma}$ is complete, $\Delta$ contains $\{1,\ldots,n\}$ and $\emptyset$. This corresponds to the cases of Remark \ref{extreme}.
\end{remark}

\begin{proposition}\label{Htrivial}
Let $\mathcal{L}$ be a line bundle on $\mathbb{P}_{\mathbf{\Sigma}}$. Then $\mathcal{L}$ is $\mathrm{H}-$trivial if and only if $\mathcal{L}$ does not lie in $FS_I$ for any $I\in \Delta$.
\end{proposition}
\begin{proof}
By Proposition \ref{Coho}, line bundle $\mathcal{L}$ is $\mathrm{H}-$trivial if and only if $$\mathrm{H}^{red}_{rk(N)-j-1}(Supp(\mathbf{r}))=0$$ for all
 $j\in \mathbb{Z}_{\geq 0}$ and all $\mathbf{r}=(r_i)_{i=1}^{n}\in \mathbb{Z}^n$ such that $\mathcal{O}(\sum_{i=1}^{n}r_iE_i)\cong \mathcal{L}$. By Definition \ref{deforbiddenset}, this is further equivalent to $\mathcal{L}\notin FS_I$ for any $I\in \Delta$.
\end{proof}

We introduce $\mathrm{Pic}_{\mathbb{R}}(\mathbb{P}_{\mathbf{\Sigma}})=\mathrm{Pic}(\mathbb{P}_{\mathbf{\Sigma}})\otimes \mathbb{R}$ which can be regarded as a quotient of $\mathbb{R}^n$ with basis elements given by $\bar E_i$ by the space of  
$\sum_{i=1}^n (m\cdot v_i )\,\bar E_i$ for $m\in M^*_{\mathbb R}$.
\begin{definition}\label{ZI}
For each $I\in \Delta$, we define the forbidden point by $$q_I=-\sum_{i\notin I}\bar E_i \in \mathrm{Pic}_{\mathbb{R}}(\mathbb{P}_{\mathbf{\Sigma}}).$$
We define a cone associated to $I$ with vertex at the origin to be $$Z_I=\sum_{i\in I}\mathbb{R}_{\geq 0}\bar E_i-\sum_{i\notin I}\mathbb{R}_{\geq 0} \bar E_i.$$
We define the forbidden cone $FC_I\subseteq \mathrm{Pic}_{\mathbb{R}}(\mathbb{P}_{\mathbf{\Sigma}}) $ by $$FC_I=q_I+Z_I.$$
\end{definition}
\begin{remark}
By definition, for any  $I\in \Delta$ the image of  $FS_I$ under the map $\mathrm{Pic}(\mathbb{P}_{\mathbf{\Sigma}})\to \mathrm{Pic}_{\mathbb{R}}(\mathbb{P}_{\mathbf{\Sigma}})$ is a subset of $FC_I$.
\end{remark}

In dimension two, the set $\Delta$ is especially simple.
\begin{example}
Let $\Sigma$ be a complete simplicial fan $\Sigma$ in $N$ with $n$ one-dimensional cones and $n$ lattice points $\{v_i\}_{i=1}^n$ chosen in each of the one-dimensional cones of $\Sigma$. In the case that $N=\mathbb{Z}^2$, the maximal cones of $\Sigma$ are $\mathbb{R}_{\geq 0}v_1+\mathbb{R}_{\geq 0}v_2, \mathbb{R}_{\geq 0}v_2+\mathbb{R}_{\geq 0}v_3,\ldots,\mathbb{R}_{\geq 0}v_n+\mathbb{R}_{\geq 0}v_1$, see Figure \ref{fig:1}. We describe $\Delta=\{\emptyset, \{1,\ldots,n\}\}\cup \{I\subset \{1,\ldots,n\}| C_I\text{ is disconnected}\}$. For example, we have $\{1,3\}\in \Delta$ if $n>3$, $\{n,2,3\}\in \Delta$ if $n>4$, but $\{1,2\}\notin \Delta$, $\{n,1,2\}\notin \Delta$ for all $n>2$, see Figure \ref{fig:1}.
\begin{figure}[H]
  \includegraphics[width=0.25\textwidth]{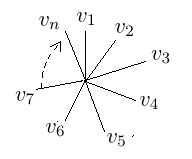}
  \caption{}
  \label{fig:1}
\end{figure}
\end{example}

\smallskip\noindent
In dimension three we can describe $\Delta$ as follows.
\begin{example}
In the case  $N=\mathbb{Z}^3$, we have
\begin{align*}
\Delta=&\{\emptyset, \{1,\ldots,n\}\}\cup \{I\subset \{1,\ldots,n\}| C_I\text{ is disconnected}\} \\
&\cup \{I\subset \{1,\ldots,n\}| C_{\{1,\ldots,n\}\setminus I}\text{ is disconnected}\}.
\end{align*}
Indeed for a line bundle $\mathcal{L}$ on $\mathbb{P}_{\mathbf{\Sigma}}$, the cohomology $\mathrm{H}^2(\mathcal{L})$ is nontrivial iff $C_I$ is disconnected and $\mathrm{H}^1(\mathcal{L})$ is nontrivial iff $C_{\{1,\ldots,n\}\setminus I}$ is disconnected.
\end{example}

We are now ready to state the first main result of this paper. We associate to any real-valued $\mathbf{\Sigma}-$piecewise linear function $\psi$ on $N_{\mathbb{R}}$ a convex polytope $\Lambda_{\psi}$ in the character space $M_{\mathbb{R}}=M\otimes\mathbb{R}$.

\begin{definition}\label{deL}
For each maximal dimensional cone $\sigma \in \mathbf{\Sigma}(m)$, let $\psi_{\sigma}\in M_{\mathbb{R}}$ be the linear function on $N_{\mathbb{R}}$ such that $\psi_{\sigma}=\psi$ in cone $\sigma$.
We define $\Lambda_{\psi}\subset M_{\mathbb{R}}$ to be the convex hull of the set $\{\psi_{\sigma}|\sigma \in \mathbf{\Sigma}(m)\}$ in the character space $M_{\mathbb{R}}$.
\end{definition}
Our first main result whose proof is given in the next section is a combinatorial condition for toric DM stacks in any dimension to have infinitely many $\mathrm{H}-$trivial line bundles.

\begin{theorem}\label{if}
Let $\mathbb{P}_{\mathbf{\Sigma}}$ be a proper smooth dimension $m$ toric DM stack associated to a complete stacky fan $\mathbf{\Sigma}=(\Sigma, \{v_i\}_{i=1}^{n})$. If there exists a $\mathbf{\Sigma}-$piecewise linear function $\psi$ such that $0<\mathrm{dim}(\Lambda_{\psi})<m$, then
there are infinitely many $\mathrm{H}-$trivial line bundles on $\mathbb{P}_{\mathbf{\Sigma}}$.
\end{theorem}
\section{Proof of the first main result}\label{H}
In this section, we prove Theorem \ref{if}.
We start by describing a key method of constructing infinitely many $\mathrm{H}-$trivial line bundles.
\begin{proposition}\label{if1}
Let $\mathbb{P}_{\mathbf{\Sigma}}$ be a proper smooth dimension $m$ toric DM stack associated to a complete stacky fan $\mathbf{\Sigma}=(\Sigma, \{v_i\}_{i=1}^{n})$ and let
$v\in \mathbf{\Sigma}(1)$ be one of the $v_i$. Suppose
there exists a globally non-linear  $\mathbf{\Sigma}-$piecewise linear function $\psi:N_{\mathbb R}\to {\mathbb R}$ which takes integer values on $N$ such that $\psi$ is constant on all lines parallel to $v$. Then
there are infinitely many $\mathrm{H}-$trivial line bundles on $\mathbb{P}_{\mathbf{\Sigma}}$.
\end{proposition}

In order to prove Proposition \ref{if1}, we need several lemmas. For an element $(a_i)_{i=1}^n $ of $\mathbb Z^n$, we define $I=\{i|a_i\geq 0\}$ and $I^c=\{1,\ldots,n\}\setminus I=\{i|a_i<0\}$.
Let $\varphi$ be a $\mathbf{\Sigma}-$piecewise linear function on $N_{\mathbb R}$ such that $\varphi(v_i)=a_i$.

\begin{lemma}\label{homotopic}
Let
$$S_{\varphi < 0}=\{v\in \big(\mathrm{N}_{\mathbb{R}}\setminus\{\mathbf{0}\}\big)\diagup \mathbb{R}_{> 0}| \varphi(v)<0\}$$ be the subset of the sphere $\big(\mathrm{N}_{\mathbb{R}}\setminus\{\mathbf{0}\}
\big)\diagup \mathbb{R}_{> 0}$ on which $\varphi$ is negative.
Let $I^c=\{i|\varphi(v_i)< 0\}$ be the subset of indices of the ray generators $v_i$ on which $\varphi$ is negative and let $C_{I^c}$ be the corresponding simplicial complex from Definition \ref{deforbiddenset}. Then $C_{I^c}$ is homotopic to $S_{\varphi < 0}$. Here we consider $C_{I^c}$ to be the geometric realization of an abstract simplicial complex \cite{A}.
\end{lemma}

\begin{proof}
We think of $C_{I^c}$ as a subspace of topological space $S_{\varphi<0}$ with the inclusion $C_{I^c}\hookrightarrow S_{\varphi<0}$ given by
\begin{equation*}
\sum_{\substack{i\in J \in C_{I^c},\\ \sum\lambda_i=1}}\lambda_i v_i\mapsto \mathbb{\mathbb{R}}_{> 0}\big(\sum_{\substack{i\in J \in C_{I^c},\\ \sum\lambda_i=1}}\lambda_i v_i\big).
\end{equation*}
Now we consider a map $F: S_{\varphi<0} \times [0,1] \rightarrow S_{\varphi<0}$ which is defined as follows. For any point $a \in S_{\varphi<0}$, let
$\sigma$ be a cone of $\mathbf{\Sigma}$ that contains $a$. We can write $a=\sum_{i\in I_{\sigma}}\lambda_i v_i $, where $I_{\sigma}=\{i|v_i \in \sigma\}$ and all $\lambda_i\geq 0$. Then we define
\begin{equation*}
F(a,t)=\sum_{\substack{i\in I_{\sigma}\\ \varphi(v_i)<0}} \lambda_i v_i + \sum_{\substack{i\in I_{\sigma}\\ \varphi(v_i)\geq 0}}(1-t)\lambda_i v_i.
\end{equation*}
Assume we choose another cone $\sigma'$ containing $a$ and write $a=\sum_{i\in I_{\sigma'}}\lambda'_i v_i $, where $I_{\sigma'}=\{i|v_i \in \sigma'\}$ and all $\lambda'_i\geq 0$. Then $\lambda_i=0=\lambda_i'$ if $v_i \notin \sigma \cap \sigma'$ and $\lambda_i=\lambda_i'$ if $v_i \in \sigma \cap \sigma'$. Thus the map $F$ is well defined. Crucially, $F$ is continuous since $\lambda_i$ change continuously when the point $a$ moves from one cone to another.
We immediately see that $F(a,0)=a$, $F(a,1) \in C_{I^c}$ for any $a\in S_{\varphi<0}$ by definition of $F$. Moreover $F(a,t)=a$ for any $a\in C_{I^c}$ and any $t\in [0,1]$ since $\varphi(v_i)<0$ for all $i\in I_{\sigma}$ if $a\in C_{I^c}$. Therefore $F$ is
a strong deformation retraction of topological space $S_{\varphi<0}$ onto subspace $C_{I^c}$.
\end{proof}

We have the following corollary.

\begin{corollary}\label{contractible}
The topological spaces $\{v\in\mathrm{N}_{\mathbb{R}}| \varphi(v)<0\}$ and $C_{I^c}$ are homotopic.
\end{corollary}

\begin{proof}
Since $\{v\in\mathrm{N}_{\mathbb{R}}| \varphi(v)<0\}=\{v\in \mathrm{N}_{\mathbb{R}}\setminus\{\mathbf{0}\}| \varphi(v)<0\}$ is a fibration with fiber $\mathbb{R}_{>0}$ over $\{v\in \big(\mathrm{N}_{\mathbb{R}}\setminus\{\mathbf{0}\}\big)\diagup \mathbb{R}_{> 0}| \varphi(v)<0\}$, these spaces are homotopic to each other. Then we use Lemma \ref{homotopic}.
\end{proof}

\begin{lemma}\label{homotopic1}
Let $S_{\varphi \geq 0}=\{v\in \big(\mathrm{N}_{\mathbb{R}}\setminus\{\mathbf{0}\}\big)\diagup \mathbb{R}_{> 0}| \varphi(v)\geq0\}$
and $I=\{i|\varphi(v_i)\geq 0\}$, then $C_{I}$ is homotopic to $S_{\varphi \geq 0}$.
\end{lemma}
\begin{proof}
The proof is analogous to that of Lemma \ref{homotopic} and is left to the reader.
\end{proof}

\begin{lemma}\label{c}
Let $I\in\{1,2,\ldots,n\}$. Then we have
$$H^{red}_{j-1}(C_{I})=\big(H^{red}_{m-j-1}(C_{I^c})\big)^*.$$
This implies that $C_I$ has nontrivial reduced homology if and only if $C_{I^c}$ has nontrivial reduced homology.
\end{lemma}
\begin{proof}
Since the sphere $S^{m-1}$ is homeomorphic to $\big(\mathrm{N}_{\mathbb{R}}\setminus\{\mathbf{0}\}\big)\diagup \mathbb{R}_{> 0}$, we have $S^{m-1}=S_{\varphi \geq 0}\sqcup S_{\varphi < 0}$. By Alexander duality (Corollary 3.45 in \cite{A}), we have an isomorphism of reduced homology and reduced cohomology $\mathrm{H}_{j-1}^{red}(S_{\varphi \geq 0})\cong \mathrm{H}_{red}^{m-1-j}(S_{\varphi < 0})$. Using the Universal Coefficient Theorem, we get $\mathrm{H}_{red}^{m-1-j}(S_{\varphi < 0})=\big(\mathrm{H}_{m-1-j}^{red}(S_{\varphi < 0})\big)^*$. Since $C_{I}$ is homotopic to $S_{\varphi \geq 0}$ by Lemma \ref{homotopic1} and $C_{I^c}$ is homotopic to $S_{\varphi < 0}$ by Lemma \ref{homotopic}, we obtain $H^{red}_{j-1}(C_{I})=\big(H^{red}_{m-j-1}(C_{I^c})\big)^*$.
\end{proof}

We still need to prove a couple of easy statements before we proceed to prove Proposition \ref{if1}.
For some fixed $v\in \mathbf{\Sigma}(1)$, we consider all lines parallel to $v$. The parametric equation of such a line is $l(t)= l(0)+tv$ for some $l(0)\in N_{\mathbb{R}}$, where $t\in \mathbb{R}$.

\begin{lemma}\label{deri}
Let $\varphi$ be a $\mathbf{\Sigma}-$piecewise linear function on $N_{\mathbb{R}}$ and $l(t)= l(0)+tv$ be a parametric equation of a line parallel to $v$. Then for any point in the interior of the cone $\sigma\in\mathbf{\Sigma}(m)$, the derivative of the function $\varphi(l(t))$ equals $\varphi_{\sigma}\cdot v$ where $\varphi_\sigma\in \mathrm{N}_{\mathbb{R}}^*$ gives the restriction of $\varphi$ to $\sigma$.
\end{lemma}

\begin{proof}
The function $\varphi$ restricts to the linear function with gradient $\varphi_\sigma$ and the claim follows.
\end{proof}

\begin{corollary}\label{conpara}
For a nonzero $v\in N_{\mathbb{R}}$, we have $\psi_{\sigma} \cdot v=0$ for all cones $\sigma \in \mathbf{\Sigma}(m)$ if and only if $\psi$ is constant on all lines parallel to $v$.
\end{corollary}

\begin{proof}[{\bf Proof of Proposition \ref{if1}.}]
Without loss of generality, we assume $v=v_1$. We claim that $\mathcal{L}=\mathcal{O}(\sum_{i=2}^{n}\psi(v_i)E_i-E_1)$ is $\mathrm{H}-$trivial. We will prove it by looking at the behavior of the relevant piecewise linear functions on lines parallel to $v_1$.

\smallskip
Let $a_i$ be the coefficient of $E_i$ in $\sum_{i=2}^{n}\psi(v_i)E_i-E_1$.
By Remark \ref{coho2},
we have $$\mathrm{H}^j(\mathbb{P}_{\mathbf{\Sigma}},\mathcal{L})=\bigoplus_{f\in M}\mathrm{H}^{red}_{m-j-1}( Supp(\mathbf{r}_f)),$$where $\mathbf{r}_f=(a_i+f \cdot v_i)_{i=1}^{n}$. In order to show $\mathrm{H}^*(\mathbb{P}_{\mathbf{\Sigma}},\mathcal{L})=0$, it is sufficient to show that $Supp(\mathbf{r}_f)$ is contractible for each $f\in M$.
We have
$$r_i:=(\mathbf{r}_f)_i= \psi(v_i) + f\cdot v_i {\rm~for~} i=2,\ldots,n{\rm~ and~} r_1=f \cdot v_1-1.$$
Let $I=\{i|r_i\geq 0\}$ and $I^c=\{1,\ldots,n\}\setminus I$.

\smallskip 
We consider two cases,  $r_1\neq-1$ and  $r_1=-1$.

\smallskip
{\bf Case $r_1\neq-1$}, equivalently  $f \cdot v_1\neq0$. We regard $\geq0$ and $<0$ as different signs. Since $f \cdot v_1$ is an integer, $f\cdot v_1 -1$ is non-negative when $f\cdot v_1\geq 1$ and $f \cdot v_1-1$ is negative when $f \cdot v_1<0$. Thus the sign of $r_i$ is the same as the sign of $f\cdot v_i+\psi(v_i)$ for all $i\in\{1,2,\ldots,n\}$. Let $\varphi$ be the $\mathbf{\Sigma}-$piecewise linear function $\varphi$ such that $\varphi(v_i)=(\psi+f)(v_i)$ for all $i\in\{1,2,\ldots,n\}$. For any line $l(t)=l(0)+tv_1$ parallel to $v_1$, we have $\varphi(l(t))=f(l(0)+tv_1)+\psi(l(t))=f(v_1)t+c_l$ for some constant value $c_l$ since $\psi$ is constant on all lines parallel to $v_1$. This implies that if $f \cdot v_1<0$ then for any line $l(t)$ parallel to $v_1$
there exists a unique $\hat t\in \mathbb R$ (which depends continuously on the line) such that $\varphi(l(t))<0$  iff $t>\hat t$. Similarly, if $f(v_1)>0$, $\varphi(l(t))$ is negative precisely when $t<\hat t$ for some  $\hat t$ that depends on $l$. Therefore $\{v\in \mathrm{N}_{\mathbb{R}}| \varphi(v)<0\}$ is contractible. So $C_{I^c}$ is contractible by Lemma \ref{contractible}. Thus $Supp(\mathbf{r}_f)=C_{I}$ has trivial reduced homology by Lemma \ref{c}.

\smallskip
{\bf Case $r_1=-1$}. We have $f(v_1)=0$. Let $\overline{\varphi}$ be the $\mathbf{\Sigma}-$piecewise linear function such that corresponds to $\mathbf{r}_f$, i.e. $\overline{\varphi}(v_i)=(\psi+f)(v_i)$ for $i\in\{2,\ldots,n\}$ and $\overline{\varphi}(v_1)=-1$. Let $l(t)=l(0)+tv_1$ be any line parallel to $v_1$. For each $l$ there exists $t_0$ such that $l(t)$ lies in the interior of a cone $\sigma_0$ with the ray $v_1$ if and only if $t>t_0$. By Lemma \ref{deri}, the derivative of $\psi(l(t))$ equals $0$ since $\psi$ is constant on all lines parallel to $v_1$. For any point in the interior of the region corresponding to $\sigma\neq\sigma_0$, the derivative of the function $\overline{\varphi}(l(t))$ equals $\overline{\varphi}_{\sigma}(v_1)=\psi_{\sigma}(v_1)+f(v_1)=\psi_{\sigma}(v_1)=0$. For any point in the interior of the region corresponding to $\sigma_0$, the derivative of the function $\overline{\varphi}(l(t))$ equals $\overline{\varphi}_{\sigma_0}(v_1)=-1$. Thus $\varphi$ is constant on $l(t)$ for $t\leq t_0$ and the derivative of the function at $l(t)$ is $(-1)$ for $t>t_0$.
So $\overline{\varphi}(l(t))$ is negative when $t>t_0$ is sufficiently large, or for all $t$. Thus $\{v\in \mathrm{N}_{\mathbb{R}}| \overline{\varphi}(v)<0\}$ is again contractible.
Then $C_{I^c}$ is contractible by Lemma \ref{contractible}. Thus $Supp(\mathbf{r}_f)=C_{I}$ has trivial reduced homology by Lemma \ref{c}.

\smallskip
Since the same argument applies to $r\psi$ instead of $\psi$, for any $r\in \mathbb{Z}$ we have
$\mathcal{L}=\mathcal{O}(\sum_{i=2}^{n}r\psi(v_i)E_i-E_1)$ is $\mathrm{H}-$trivial. These line bundles are non-isomorphic to one another because $\psi$ is non-linear.
Thus there are infinitely many $\mathrm{H}-$trivial line bundles on $\mathbb{P}_{\mathbf{\Sigma}}$.
\end{proof}

In order to prove Theorem \ref{if},  we will show that $\mathrm{dim}\Lambda_{\psi}<m$ implies that there exists $v\in \mathbf{\Sigma}(1)$ such that $\Lambda_{\psi} \cdot v=0$.

\begin{proposition}\label{Lambda} The inequality
$\mathrm{dim}\Lambda_{\psi}<m$ holds  if and only if there exist $w\in \mathbf{\Sigma}(1)$  and a constant $c\in {\mathbb R}$ such that $\Lambda_{\psi} \cdot w= c$.
\end{proposition}

\begin{proof}
We have that $\mathrm{dim}\Lambda_{\psi}<m$ if and only if $\Lambda_{\psi}$ is inside an affine hyperplane in $M_{\mathbb{R}}$. Equivalently there exists a non-zero $v\in N_{\mathbb{R}}$ such that $\langle\Lambda_{\psi},v\rangle=c$, where $c$ is a constant. The essence of this proposition is that $v$ may be chosen in $\mathbf{\Sigma}(1)$. The ``if'' part is clear, so we focus on the ``only if'' part.

\smallskip
As a first step, we adjust $\psi$ by adding a linear function in $M$ such that $\psi_{\sigma}(v)=0$ for all cones $\sigma \in \mathbf{\Sigma}(m)$ and $\psi=0$ on the (unique) cone $\tau$
which contains $v$ in its interior. This clearly has no effect on the statement of the proposition, since adding a global linear function amounts to a parallel translation of $\Lambda_{\psi}$. Another consequence of this shift is that $\psi$ is now constant on lines parallel to $v$.

\smallskip
We now claim that for every $\sigma \in \mathbf{\Sigma}(m)$, there exists a
$\widetilde{\sigma}\in \mathrm{Star}(\tau)(m):=\{\sigma \in \mathbf{\Sigma}(m)| \tau\subseteq \sigma\}$ such that $\psi_{\sigma}=\psi_{\widetilde{\sigma}}$. To prove it, consider the projection $pr: N_{\mathbb{R}}\rightarrow N_{\mathbb{R}}/ \mathbb{R}v$ and choose a maximum-dimensional cone $\widetilde{\sigma}\in \mathrm{Star}(\tau)(m)$ such that $pr(\sigma)$ and $pr(\widetilde{\sigma})$ have overlapping interiors. Let $D=pr(\sigma)\cap pr(\widetilde{\sigma})$, see Figure \ref{fig:2}.
There are linear functions $\varphi_{pr(\sigma)}$ and $\varphi_{pr(\widetilde{\sigma})}$ on $N_{\mathbb{R}}/ \mathbb{R}v$ such that $\psi_{\sigma}=\varphi_{pr(\sigma)}\circ pr$ and $\psi_{\widetilde{\sigma}}=\varphi_{pr(\widetilde{\sigma})}\circ pr$.

\smallskip
For any $p\in D$, we pick a point $q_1$ in $\sigma$ such that $pr(q_1)=p$ and a point $q_2$ in $\widetilde{\sigma}$ such that $pr(q_2)=p$. Since $q_1$ and $q_2$ are on a line parallel to $v$, by Lemma \ref{conpara}, we get $\psi(q_1)=\psi(q_2)$. Since $\psi_{\sigma}(q_1)=\psi(q_1)$ and $\psi_{\widetilde{\sigma}}(q_2)=\psi(q_2)$, we have $$\varphi_{pr(\sigma)}(p)=\psi_{\sigma}(q_1)=\psi_{\widetilde{\sigma}}(q_2)=\varphi_{pr(\widetilde{\sigma})}(p).$$ So we get $\varphi_{pr(\sigma)}=\varphi_{pr(\widetilde{\sigma})}$ on $D$ which is a full-dimensional set. This implies $\varphi_{pr(\sigma)}=\varphi_{pr(\widetilde{\sigma})}$ on $N_{\mathbb{R}}/ \mathbb{R}v$. Thus $\psi_{\sigma}=\psi_{\widetilde{\sigma}}$.
\begin{figure}[H]
  \includegraphics[width=0.28\textwidth]{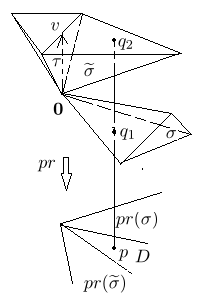}
  \caption{}
  \label{fig:2}
\end{figure}
 We choose a vertex $w$ of the cone $\tau$. For any maximum-dimensional cone $\widetilde{\sigma}$ in $\mathrm{Star}(\tau)$, we have $\psi_{\widetilde{\sigma}}(w)=\psi(w)=0$ since $w\in \tau \subset \widetilde{\sigma}$. Then for all cones $\sigma \in \mathbf{\Sigma}(m)$, we obtain $\psi_{\sigma}(w)=\psi_{\widetilde{\sigma}}(w)=0$.
\end{proof}

Now we are ready to prove Theorem \ref{if}.

\begin{proof}[{\bf Proof of Theorem \ref{if}.}]
By Proposition \ref{Lambda}, there exists $w\in \mathbf{\Sigma}(1)$ such that 
 $\Lambda_{\psi} \cdot w= c$ for some constant $c$. By adding a global linear function, we may ensure that $c=0$. Corollary \ref{conpara} and Proposition \ref{if1} then imply the result.
\end{proof}

\section{Dimension three case}\label{three}
In this section we focus on proper smooth dimension three toric DM stacks $\mathbb{P}_{\mathbf{\Sigma}}$ associated to a complete stacky fan $\mathbf{\Sigma}=(\Sigma, \{v_i\}_{i=1}^{n})$ in the lattice $N$ of rank three. We give a criterion for when there exist infinitely many $\mathrm{H}-$trivial line bundles on $\mathbb{P}_{\mathbf{\Sigma}}$ for smooth toric varieties and DM stacks in dimension three under the (technical) assumption that there is no more than one pair of collinear rays in $\mathbf{\Sigma}$.

\smallskip
The arguments of this section are rather cumbersome. We start with several preliminary results.
The following easy lemma highlights the importance of diagonals in $\mathbf{\Sigma}$, i.e., pairs $(s,t)$ such that $v_s$ and $v_t$ are collinear.

\begin{lemma}\label{g}
Let $\bar E=\sum_{i=1}^{n}f(v_i)
\bar E_i$ be an element in $\mathrm{Pic}_{\mathbb{R}}(\mathbb{P}_{\mathbf{\Sigma}})$, where $f$ is a $\mathbf{\Sigma}-$piecewise linear function on $N_{\mathbb{R}}$. Let $s\in \{1,2,\ldots,n\}$. Then there exists a $\mathbf{\Sigma}-$piecewise linear function $g$ on $N_{\mathbb{R}}$ such that
 \begin{itemize}
 \item $\bar E=\sum_{i=1}^{n}g(v_i)\bar E_i$,
 \item $g(v_s)=0$,
 \item $g(v_i)\neq 0$ for all $v_i\notin \mathbb{R}v_s$,
 \item If $v_t\in \mathbb{R}v_s$, $t\neq s$, then $g(v_t)=0$ iff $f$ is linear on $\mathbb{R}v_s$.
 \end{itemize}
 \end{lemma}
 \begin{proof}
We consider $g(v)=f(v)-m\cdot v$, where $m$ is a generic element in the affine plane $\{m\in M_{\mathbb{R}}| m\cdot v_s=f(v_s)\}$. We immediately get $g(v_s)=0$ and  $\bar E=\sum_{i=1}^{n}g(v_i)\bar E_i$.
 The fact that $m$ is  generic means that $m\cdot v_t\neq f(v_t)$ for all $i$ except possibly when $v_t\in \mathbb{R}v_s$. If $v_t\in \mathbb{R}v_s$, then $g$ is linear on $\mathbb{R}v_s$ iff $g(v_t)=0$, and g is linear on $\mathbb{R}v_s$ iff $f$ is linear on $\mathbb{R}v_s$.
 \end{proof}

\begin{definition}
For a point $v_i\in\{v_i\}_{i=1}^{n}$, we denote the neighborhood of $v_i$ to be $B_i=\{v_j\in\{v_1,\ldots,v_n\}|v_j \text{ and } v_{i} \text{ span a two-dimensional cone of } \mathbf{\Sigma} \}.$
\end{definition}

\begin{definition}
Let $v_{j_1}, v_{j_2}, \ldots, v_{j_l}$ be the vectors in $B_i$ which are ordered clockwise when looking from $v_i$ to the origin. Let $g$ be a $\mathbf{\Sigma}-$piecewise linear function on $N_{\mathbb{R}}$. We regard $\geq0$ and $<0$ as different signs. We count the number of pairs of vectors $\{v_{j_{k}},v_{j_{k+1}}\}\subset \{v_{j_1}, \ldots, v_{j_l}\}$ such that $f(v_{j_{k}})$ and $f(v_{j_{k+1}})$ have different signs. We call it the number of sign changes of $f$ among $v_{j_1}, \ldots, v_{j_l}$. For example, there are exactly two sign changes in $B_i$ in Figure \ref{fig:3}.
\begin{figure}[H]
  \includegraphics[width=0.32\textwidth]{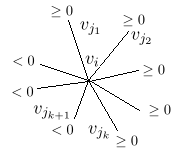}
  \caption{}
  \label{fig:3}
\end{figure}
\end{definition}

Recall that in Definition \ref{ZI}, the cones $Z_I$ are shifts of the forbidden cones $FC_I$ to the origin. The following easy observation will be used below. 
\begin{remark}\label{interiorZI}
An element $\bar E$ of  $\mathrm{Pic}_{\mathbb{R}}$ given by
$$
\bar E=\sum_i r_i \bar E_i
$$
with $r_i>0$ for all $i\in I$ and $r_i<0$ for all $i\not\in I$ lies in the interior of $Z_I$. Indeed, any positive linear combination of generators of a polyhedral cone lies in its interior.
\end{remark}

The nontrivial Lemma  \ref{two} below, which will be key to subsequent arguments, establishes some peculiar features of elements in $\mathrm{Pic}_{\mathbb R}(\mathbb{P}_{\mathbf{\Sigma}})$ which are not contained in the interior of any $Z_I$. Our interest in such elements of $\mathrm{Pic}_{\mathbb R}(\mathbb{P}_{\mathbf{\Sigma}})$ stems from the fact that if $\mathbb{P}_{\mathbf{\Sigma}}$ has infinitely many $\mathrm{H}-$trivial line bundles, then their accumulation points in the sphere $(\mathrm{Pic}(\mathbb{P}_{\mathbf{\Sigma}})\setminus\{0\})/{\mathbb R}_{>0}$ have this property (see \cite{W}).

\begin{lemma}\label{two}
Let $\bar E=\sum_{i=1}^{n}g(v_i)\bar E_i$ be a nonzero element in $\mathrm{Pic}_{\mathbb R}(\mathbb{P}_{\mathbf{\Sigma}})$ which is not contained in the interior of $Z_{I}$ for any $I\in \Delta$, where $g$ is a $\mathbf{\Sigma}-$piecewise linear function on $N_{\mathbb{R}}$. Assume that there exists some $s\in \{1,2,\ldots,n\}$ such that $g(v_{s})=0$ and $g(v_i)\neq 0$ for all $i\neq s$. Then there exist exactly two sign changes in $B_{s}$, see Figure \ref{fig:3}.
\end{lemma}
\begin{proof}
The main idea is that one can perturb $g$ by a small global linear function to achieve $g(v_s)>0$ or $g(v_s)<0$ without changing signs of the other $g(v_i)$. Such perturbation changes the relevant simplicial complexes slightly and in most cases one of them will have a non-trivial reduced homology which would be inconsistent with $\bar E$ not being in the interior of any $Z_I$.

\smallskip
Let $C_{> 0}$ be the simplicial complex $C_{\{i|g(v_i)>0\}}$ and $C_{\geq0}$ be the simplicial complex $C_{\{i|g(v_i)\geq 0\}}$. The reduced simplicial homology complex $W(C_{>0})$ is as follows:
\begin{equation*}
0\rightarrow\mathop{\bigoplus}_{\substack{J\in C_{>0} \\ |J|=3}}\mathbb{C}\rightarrow\mathop{\bigoplus}_{\substack{J\in C_{>0} \\ |J|=2}}\mathbb{C}\rightarrow\mathop{\bigoplus}_{\substack{J\in C_{>0} \\ |J|=1}}\mathbb{C}\rightarrow \mathbb{C}\rightarrow 0
\end{equation*}
The reduced simplicial homology complex $W(C_{\geq0})$ is as follows:
\begin{equation*}
0\rightarrow\mathop{\bigoplus}_{\substack{J\in C_{\geq0} \\ |J|=3}}\mathbb{C}\rightarrow\mathop{\bigoplus}_{\substack{J\in C_{\geq0} \\ |J|=2}}\mathbb{C}\rightarrow\mathop{\bigoplus}_{\substack{J\in C_{\geq0} \\ |J|=1}}\mathbb{C}\rightarrow \mathbb{C}\rightarrow 0
\end{equation*}
We use inclusions $\mathop{\bigoplus}_{\substack{J\in C_{>0}\\ |J|=k}}\mathbb{C}\hookrightarrow\mathop{\bigoplus}_{\substack{J\in C_{\geq0}\\ |J|=k}}\mathbb{C}$ for each $k$ to obtain an exact sequence of complexes:
\begin{equation}\label{exact}
0\rightarrow W(C_{>0})\rightarrow W(C_{\geq0})\rightarrow W(C_{\geq0})/ W(C_{>0})\rightarrow 0
\end{equation}

\smallskip
We now take a linear function $f$ such that $f(v_{s})>0$ and $|f(v_i)|<|g(v_i)|$ for all $i\neq s$. Let $g'=g+f$,  we have $g'(v_i)=(g+f)(v_i)\neq 0$ for $i\in \{1,2,\ldots,n\}$. Thus
$$\bar E=\sum_{i\in I}\alpha_i\bar E_i+\sum_{i\notin I}\alpha_i\bar E_i$$
in $\mathrm{Pic}_{\mathbb{R}}(\mathbb{P}_{\mathbf{\Sigma}})$, where $\alpha_i=g'(v_i)\neq 0$ for all $i$ and $I:=\{i|g'(v)>0\}$. We also have $C_{\geq0}=C_{I}$ by the definition $C_I$.
Since $\mathcal{L}$ is not contained in the interior of $Z_{I}$ for any $I\in \Delta$, we get $I\notin \Delta$, in view of Remark \ref{interiorZI}. This implies $C_{\geq0}=C_I$ has trivial reduced homology.
Similarly, we consider $g'=g-f$ to show that $C_{>0}$ has trivial reduced homology.

\smallskip
Since $W(C_{\geq0})$ and $W(C_{>0})$ have trivial homology, by snake lemma and the exact sequence \eqref{exact}, the complex $W(C_{\geq0})/ W(C_{>0})$ has trivial homology. Thus the Euler characteristic $\chi(W(C_{\geq0})/ W(C_{>0}))=0$.

\smallskip
If $g(v_i)<0$ for all $v_i\in B_{s}$, then 
$\chi(W(C_{\geq0})/ W(C_{>0}))=1$, because the only nonzero space in $W(C_{\geq0})/ W(C_{>0}))$ is the contribution of the zero-dimensional cell $\{s\}$. If $g(v_i)>0$ for all $v_i\in B_{s}$, then we also have 
 $\chi(W(C_{\geq0})/ W(C_{>0}))= 1$ by observing that there is $1$ zero-dimensional cell, $|B_s|$ one-dimensional cells and $|B_s|$ two-dimensional cells.

\smallskip
Thus there exist at least two sign changes in $B_{s}$. A length $l$ connected component of positive sign in $B_s$ provides $l$ one-dimensional cells and $(l-1)$ two-dimensional cells, so the Euler characteristic $\chi(W(C_{\geq0})/ W(C_{>0}))$ equals
$$1-\text{ the number of connect components of positive sign in } B_s .$$
Thus there exists only one component of positive sign in $B_s$. This implies that there are exactly two sign changes in $B_s$.
\end{proof}

Let $\bar E=\sum_{i=1}^{n}g(v_i)\bar E_i$ be an element in $\mathrm{Pic}_{\mathbb R}(\mathbb{P}_{\mathbf{\Sigma}})$ which is not contained in the interior of $Z_{I}$ for any $I\in \Delta$, where $g$ is a $\mathbf{\Sigma}-$piecewise linear function on $N_{\mathbb{R}}$. Assume there is some $s\in \{1,2,\ldots,n\}$ such that $g(v_{s})=0$. Then we can consider the projection $\pi:N_{\mathbb{R}}\rightarrow\mathbb{R}^2=N_{\mathbb{R}}/\langle v_{s}\rangle$ along the line $\mathbb{R}v_{s}$.

\smallskip
Under the projection, the neighborhood $B_{s}$ gives a two-dimensional complete stacky fan $\mathbf{\Sigma}'=(\Sigma', \{w_i\}_{i\in J})$, where $w_i=\pi(v_i)$, $v_i\in B_s$ and $J=\{i\in\{1,2,\ldots,n\}|v_i\in B_s\}$.
We can think of $g$ as a $\mathbf{\Sigma}'-$piecewise linear function on $\mathbb{R}^2=N_{\mathbb{R}}/\langle v_{s}\rangle$ since $g(v_{s})=0$. Let $\{\bar E'_i\}_{i\in J}$ be the generators of $\mathrm{Pic}_{\mathbb R}(\mathbb{P}_{\mathbf{\Sigma}'})$. Let $\Delta'=\{\emptyset, J\}\cup \{I'\subset J | C'_{I'}\text{ is disconnected}\}$. For any $I'\in\Delta'$, let $Z'_{I'}=\sum_{i\in I'}\mathbb{R}_{\geq 0}\bar E'_i-\sum_{i\notin I'}\mathbb{R}_{\geq 0} \bar E'_i$ be the cone associated to $I'$ with vertex at the origin which is obtained by shifting the forbidden cone $FC'_I\subseteq\mathrm{Pic}_{\mathbb R}(\mathbb{P}_{\mathbf{\Sigma}'})$.
\begin{lemma}\label{doto2}
The point $\bar D=\sum_{i\in J}g(v_i)\bar E'_i\in \mathrm{Pic}(\mathbb{P}_{\mathbf{\Sigma'}})$ is not contained in the interior of $Z'_{I'}$ for any $I'\in\Delta'$.
\end{lemma}
\begin{proof}
Assume $\bar D$ is contained in $Z'_{I'}$ for some $I'\in\Delta'$. This implies that there exists a linear function $f$ on $\mathbb{R}^2=\mathbb{R}^3/\langle v_{s}\rangle$ such that $f(w_i)+g(v_i)\neq 0$ for any $i\in J$ and $I'=\{i\in J | f(w_i)+g(v_i) \}$ has non-trivial reduced homology. This would mean that the number of sign changes in $B_s$ is not two, which contradicts Lemma \ref{two}.
\end{proof}

\begin{lemma}\label{half}
Let $\bar E=\sum_{i=1}^{n}g(v_i){\bar E}_i$ be an element in $\mathrm{Pic}(\mathbb{P}_{\mathbf{\Sigma}})$ which is not contained in the interior of $Z_{I}$ for any $I\in \Delta$, where $g$ is a $\mathbf{\Sigma}-$piecewise linear function on $N_{\mathbb{R}}$. Assume $v_s$ is such that there is no $v_t$ with $t \neq s$, $\mathbb{R}v_s=\mathbb{R}v_t$; or $g$ is not linear on $\mathbb{R}v_s$. Then
\begin{itemize}
\item either $g$ is linear on $\mathrm{Star}(v_s)=\{\sigma \in \mathbf{\Sigma}| v_s\subseteq \sigma\}$
\item or $g|_{\mathrm{Star}(v_s)}$ has two half spaces of linearity, i.e., there exists $v_p, v_q\in B_s$ such that $v_p, v_s, v_q, \mathbf{0}$ are coplanar and $g|_{\mathrm{Star}(v_s)}$ is linear on either side of this plane, see Figure \ref{fig:4}.
\end{itemize}
\begin{figure}[H]
  \includegraphics[width=0.23\textwidth]{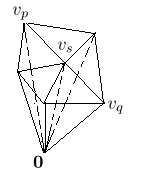}
  \caption{}
  \label{fig:4}
\end{figure}
\end{lemma}
\begin{proof}

By Lemma \ref{g}, we can assume $g(v_s)=0$. Then by Lemma \ref{doto2}, we have
$\bar D=\sum_{i\in J}g(v_i)\bar E'_i\in \mathrm{Pic}_{\mathbb R}(\mathbb{P}_{\mathbf{\Sigma'}})$ not in the interior of $Z'_{I'}$ for any $I\in \Delta'$.
If $\bar D=0$ in $\mathrm{Pic}_{\mathbb{R}}(\mathbb{P}_{\mathbf{\Sigma'}})$, then $g$ is linear on $\mathrm{Star}(v_s)$. Otherwise, by Lemma 4.1 in the paper \cite{W},  $g|_{\mathrm{Star}(v_s)}$ is a pullback of a piecewise linear function on $N_{\mathbb{R}}/\mathbb{R}v_{s}$ with two half plane regions of linearity.
Thus $g|_{\mathrm{Star}(v_s)}$ has two half spaces of linearity, see Figure \ref{fig:4}.
\end{proof}

We will now consider the relatively easy case of dimension three fans with no collinear rays.
\begin{proposition}\label{nodia}
Assume that $v_i$ and $v_j$ are not collinear for any $\{i, j\}\subset\{1,2,\ldots,n\}$ and there are infinitely many $\mathrm{H}-$trivial line bundles on $\mathbb{P}_{\mathbf{\Sigma}}$. Then there exists a plane which does not intersect with the interior of any maximal cone of $\mathbf{\Sigma}$.
\end{proposition}

\begin{proof}
Since there are infinitely many $\mathrm{H}-$trivial line bundles on $\mathbb{P}_{\mathbf{\Sigma}}$, by the argument of proof in Proposition 3.8 in \cite{W}, there exists a non-zero element $\bar E=\sum_{i=1}^{n}g_i\bar E_i\in \mathrm{Pic}_{\mathbb R}( \mathbb{P}_{\mathbf{\Sigma}})$ which is not contained in the interior of any $Z_{I}$ for $I\in \Delta$. \footnote{The idea of the argument  of \cite{W} is that the infinite set of $\mathrm{H}-$trivial line bundles will have an accumulation point in $\Big(\mathrm{Pic}_{\mathbb R}(\mathbb{P}_{\mathbf{\Sigma}})\setminus \{0\}\Big)/{\mathbb R}_{>0}$. If this point lies in the interior of some $Z_I$ then we get $\mathrm{H}-$trivial line bundles deep inside $Z_I$, which is impossible.
}
We will introduce a $\mathbf{\Sigma}-$piecewise linear function $g$ on $N_{\mathbb{R}}$
with $g_i=g(v_i)$.
\begin{figure}[H]
  \includegraphics[width=0.3\textwidth]{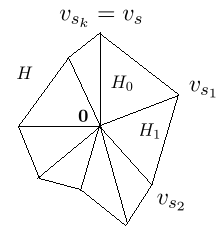}
  \caption{}
  \label{fig:5}
\end{figure}
Since $\bar E\neq 0$, there exists a vertex $v_s$ such that $g$ is not linear in $\mathrm{Star}(v_s)$. By Lemma \ref{half}, $g$ has two half planes of linearity in $\mathrm{Star}(v_s)$. The line breaking the linearity corresponds to a plane in $\mathbb{R}^3$ which we denote by $H_0$. We know that $H_0$ passes through the vector $v_s$, the origin $\mathbf{0}$ and a point in $B_s$ which we denote by $v_{s_1}$. Also, $g$ is linear on each side of plane $H_0$ in $\mathrm{Star}(v_s)$. Thus $H_0$ does not intersect with interior of any maximal cones of $\mathbf{\Sigma}$ with $v_s$ as a ray. By assumption, we know $v_{s_1}$ and $v_i$ is not collinear for any $i\neq s_1$. By Lemma \ref{g} and the same argument for $v_s$, we get a plane $H_1$ passing through $v_{s_1}$ such that $g$ is linear on each side of plane $H_1$ in $\mathrm{Star}(v_s)$. Thus $H_1$ equals $H_0$ and does not intersect with interior of any maximal cones of $\mathbf{\Sigma}$ with $v_{s_1}$ as ray. We continue the process, and eventually come around to find a positive number $k$ such that $v_{s_k}=v_s$. We get a plane $H=H_i$ for $i=0,\ldots,k$ which does not intersect with interior of any maximal cones of $\mathbf{\Sigma}$. See Figure \ref{fig:5}.
\end{proof}

\begin{theorem}\label{onlyif}
Let $\mathbb{P}_{\mathbf{\Sigma}}$ be a proper smooth dimension three toric DM stack associated to a complete stacky fan $\mathbf{\Sigma}=(\Sigma, \{v_i\}_{i=1}^{n})$. Assume there exist no collinear pairs of rays in $\mathbf{\Sigma}$ and
there are infinitely many $\mathrm{H}-$trivial line bundles on $\mathbb{P}_{\mathbf{\Sigma}}$. Then there is a $\mathbf{\Sigma}-$piecewise linear function $\psi$ which takes integer values on $N$ such that $\mathrm{dim}(\Lambda_{\psi})<3$.
\end{theorem}

\begin{proof}
By Proposition \ref{nodia}, there is a plane $H$ which does not intersect the interior of any maximal cone of $\mathbf{\Sigma}$. Then a continuous piecewise linear function $\psi$ which is zero on one side of $H$ and is nonzero on the other side of it,  has $\mathrm{dim}(\Lambda_{\psi})=1$. Since $H$ is rational, a scalar multiple of $\psi$ takes integer values on $N$.
\end{proof}


Now we consider the case when there exists exactly one collinear pair in $\{v_1,\ldots,v_n\}$.

\begin{figure}[H]
  \includegraphics[width=0.35\textwidth]{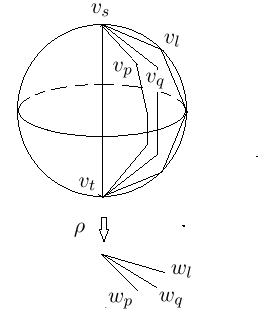}
  \caption{}
  \label{fig:6}
\end{figure}

\begin{lemma}\label{dia}
Assume that $v_s$ and $v_t$ form the only collinear pair and there are infinitely many $\mathrm{H}-$trivial line bundles on $\mathbb{P}_{\mathbf{\Sigma}}$. Then either there is a plane which does not intersect  the interior of any maximal cone of $\mathbf{\Sigma}$ or there exist at least three half planes which contain $v_s$ and $v_t$ and do not intersect with the interior of any maximal cones of $\mathbf{\Sigma}$.
\end{lemma}
\begin{proof}
Since there are infinitely many $\mathrm{H}-$trivial line bundles on $\mathbb{P}_{\mathbf{\Sigma}}$, by the argument of proof of Proposition 3.8 in \cite{W}, there exists a non-zero element $\bar E=\sum_{i=1}^{n}g(v_i)\bar E_i\in\mathrm{Pic}_{\mathbb R}( \mathbb{P}_{\mathbf{\Sigma}})$ which is not contained in any $Z_{I}$ for $I\in \Delta$, where $g$ is a $\mathbf{\Sigma}-$piecewise linear function on $\mathbb{R}^3$.

\smallskip
If $g$ is not linear along the line $\mathbb{R}v_s$, then by Lemma \ref{g}, it is the same case as the one without a collinear pair. By the same argument as in Proposition \ref{nodia}, there is a plane which does not intersect with interior of any maximal cones of $\mathbf{\Sigma}$.

\smallskip
 Now we consider the case when $g$ is linear along the line $\mathbb{R}v_s=\mathbb{R}v_t$. If $g$ has two half planes of linearity in $\mathrm{Star}(v_s)$ and $\mathrm{Star}(v_t)$, then with the same argument as in the case without collinear pairs in Proposition \ref{nodia} we deduce that there is a plane which does not intersect with the interior of any maximal cones of $\mathbf{\Sigma}$ with $v_s$ as a ray. Thus without loss of generality, we assume $g$ has at least three regions of linearity in $\mathrm{Star}(v_s)$. Thus at least three vectors in $B_s$ break the linearity. Thus in $B_s$, we can pick three $v_p$, $v_q$ and $v_l$ such that $g$ is given by a linear function in the cone spanned by $\{v_s,v_p,v_q\}$ and is given by another linear function in the cone span by $\{v_s,v_q,v_l\}$. This implies the two-dimensional cone $C_{sq}$ spanned by $\{v_s, v_q\}$ break the linearity.

\smallskip
Now we consider $v_q$ and its neighborhood. By assumption, we know $v_q$ and $v_i$ is not collinear for any $i\neq p$. By Lemma \ref{g} and the same argument as in the proof of Proposition \ref{nodia}, we get a plane $H_0$ passing through $v_{q}$ such that $g$ is linear on each side of plane $H_0$ near $v_q$. Thus $H_0$ is the plane of the cone $C_{sp}$ spanned by $\{v_s, v_q\}$ and does not intersect with interior of any maximal cones of $\mathbf{\Sigma}$ with $v_q$ as ray. We continue the process, there exits a positive number $k$ such that $v_{q_k}=v_t$. We get a half plane $H_q$ which does not intersect with interior of any maximal cones of $\mathbf{\Sigma}$ and passes through $v_s$.

\smallskip
Similarly, we get half planes $H_p$ and $H_l$ which do not intersect with interior of any maximal cones of $\mathbf{\Sigma}$ and passes through $v_s$. See Figure \ref{fig:6}.
\end{proof}

\begin{theorem}\label{oneif}
Let $\mathbb{P}_{\mathbf{\Sigma}}$ be a proper smooth dimension three toric DM stack associated to a complete stacky fan $\mathbf{\Sigma}=(\Sigma, \{v_i\}_{i=1}^{n})$. Assume there exists only one collinear pair of rays in $\mathbf{\Sigma}$ and
there are infinitely many $\mathrm{H}-$trivial line bundles on $\mathbb{P}_{\mathbf{\Sigma}}$. Then there is a $\mathbf{\Sigma}-$piecewise linear function $\psi$ such that $\mathrm{dim}(\Lambda_{\psi})<3$.
\end{theorem}
\begin{proof}
Let $v_s$ and $v_t$ be the only collinear pair.
By assumption and Proposition \ref{dia}, either there exists a projection $\pi: \mathbb{R}^3\rightarrow \mathbb{R}^1$ along a plane which does not intersect with the interior of any maximal cones of $\mathbf{\Sigma}$ or there are three vectors $v_p,v_q,v_l \in B_s$ and a projection $\rho: \mathbb{R}^3\rightarrow \mathbb{R}^2$ alone the line $\mathbb{R}v_s$ such the $H_i=\rho^{-1}\rho(v_i)$ is a half plane which does not intersect with interior of any maximal cones of $\mathbf{\Sigma}$ for each $i=p,q,l$. In former case, with the same argument in Theorem \ref{onlyif}, there is a $\mathbf{\Sigma}-$piecewise linear function which is constant on all lines parallel to $\mathbb{R}v_s$.

\smallskip
Let $w_i=\rho(v_i)$ for each $i=p,q,l$.
In latter case, we pick a function $f$ on $\mathbb{R}^2$ which is linear respectively on the cone $C_{pq}$ spanned by $\{w_p,w_q\}$, the cone $C_{ql}$ spanned by $\{w_q,w_l\}$ and the cone $C_{lp} $ spanned by $\{w_l,w_p\}$. Let $\psi=\rho^*(f)$, see Figure \ref{fig:6}. We get $\psi=\rho^*(f)$ is a $\mathbf{\Sigma}-$piecewise linear function on $\mathbb{R}^3$ since the preimage of $\rho(v_i)$ is $H_i$ which does not intersect with interior of any maximal cones of $\mathbf{\Sigma}$ for each $i=p,q,l$. Also, $\psi$ is a $\mathbf{\Sigma}-$piecewise linear function which is constant on any line parallel to $v_s$ since the image of any line under $\rho$ is a point.
Then Corollary \ref{conpara} and Proposition \ref{Lambda} imply the result.
\end{proof}

Putting all of it together, we get our second main result.
\begin{theorem}\label{oneonly}
Let $\mathbb{P}_{\mathbf{\Sigma}}$ be a proper smooth dimension three toric DM stack associated to a complete stacky fan $\mathbf{\Sigma}=(\Sigma, \{v_i\}_{i=1}^{n})$. Assume there exists no more than one collinear pair of rays in $\mathbf{\Sigma}$. Then
there are infinitely many $\mathrm{H}-$trivial line bundles on $\mathbb{P}_{\mathbf{\Sigma}}$ if and only if there exists a $\mathbf{\Sigma}-$piecewise linear function $\psi$ such that $0<\mathrm{dim}(\Lambda_{\psi})<3$.
\end{theorem}
\begin{proof}
Theorem \ref{onlyif}, Theorem \ref{oneif} and Theorem \ref{if} imply the result.
\end{proof}

\section{Comments}\label{comments}
In this section, we discuss possible refinements and generalizations of the results of this paper. The dimension of $\mathbb{P}_{\mathbf{\Sigma}}$ is now arbitrary.

\smallskip
For a proper smooth dimension $m$ toric DM stack $\mathbb{P}_{\mathbf{\Sigma}}$, we have the following three statements:
\begin{enumerate}
\item There exists a $\mathbf{\Sigma}-$piecewise linear function $\psi:N_{\mathbb R} \to {\mathbb R}$ such that $0<\mathrm{dim}(\Lambda_{\psi})<m$.
\item There are infinitely many $\mathrm{H}-$trivial line bundles on $\mathbb{P}_{\mathbf{\Sigma}}$.
\item There is a nonzero element $\bar E \in \mathrm{Pic}_{\mathbb R}(\mathbb{P}_{\mathbf{\Sigma}})$ which is not contained in the interior of $Z_{I}$ for any $I\in \Delta$.
\end{enumerate}
By Theorem \ref{if}, $(1)$ implies $(2)$ without any additional assumptions. 
By the argument of \cite{W}, $(2)$ also implies $(3)$ in full generality. By Theorem \ref{oneonly}, $(2)$ implies $(1)$ when $m=3$ and there exist no more than one collinear pair of rays in $\mathbf{\Sigma}$. By \cite{W} all three conditions are equivalent for $m=2$.

\smallskip
In the most optimistic scenario, $(3)$ implies $(1)$ in full generality, so all three conditions above are equivalent, but we currently have no methods to approach this problem, even in dimension three. Alternatively, it would be also interesting to find an example of a proper toric DM stack that satisfies  $(3)$ but not $(1)$.

\smallskip
Naturally, other questions to ask is whether $(2)$ implies $(1)$ in full generality or if $(3)$ implies $(2)$. At present, we can neither prove these implications nor provide counterexamples. Still, we hope that the methods and the approach of this paper will be useful for settling these mysteries.

\end{document}